\documentclass[12pt]{article}
\usepackage{amsmath,amsfonts,amssymb,amsthm}
\usepackage{natbib}

\def\spacingset#1{\renewcommand{\baselinestretch}%
{#1}\small\normalsize} \spacingset{1}

\RequirePackage[colorlinks,citecolor=blue,urlcolor=blue]{hyperref}

\usepackage{multicol} 
\usepackage{dcolumn}
\usepackage{bm}
\usepackage{multirow}
\usepackage{lscape}
\usepackage{setspace}
\usepackage{rotating}
\usepackage{changepage}
\usepackage[latin1]{inputenc}
\usepackage[titletoc]{appendix}

\newcommand{\Pmeas}{\mathbb{P}}

\newcommand{\brac}[1]{\left ( #1 \right )}
\newcommand{\brak}[1]{\left [ #1 \right ]}
\newcommand{\brat}[1]{\left \{ #1 \right \}}

\newcommand{\norm}[1]{\left\Vert#1\right\Vert}

\newcommand{\abs}[1]{\lvert #1 \rvert}

\newcommand{\vc}{\text{vec}}

 \newcommand{\T}{\top}
\newcommand{\reals}{\mathbb{R}}
\usepackage{caption,subcaption}
\newcolumntype {d}[1]{D{.}{.}{#1}} \newcolumntype {.}{D{.}{.}{-1}}

\usepackage{color}
\definecolor {orange}{rgb}{1,0.5,0}

\newtheorem{theorem}{Theorem}[section]

\newtheorem{corollary}{Corollary}
\newtheorem{proposition}[theorem]{Proposition}
\theoremstyle{remark}
\newtheorem{definition}[theorem]{Definition}
\newtheorem{example}{Example}

\newtheorem{assumption}[theorem]{Assumption}

\DeclareMathOperator{\tr}{tr}
\DeclareMathOperator{\vecme}{vec}

\newcommand{\PP}{\mathcal P}
\newcommand{\QQ}{\mathcal Q}
\newcommand{\RR}{\mathcal R}
  
\begin{document}

\title{Constrained polynomial likelihood\thanks{We are grateful to Don Andrews, Oleg Bondarenko, Albina Danilova, Marc Van Uffelen, and participants at the 2021 SoFiE Conference (San Diego), 2021 Brazilian Meeting of Finance, Lunchtime Workshop at Boston University, and the Remote Seminar Series on Computational Economics and Finance for helpful comments. Paul Schneider gratefully acknowledges the SNF grant 
100018\_189086 ``Scenarios".}}

\author{Caio Almeida\thanks{calmeida@princeton.edu, 
      Department of Economics and BCF, 
      Princeton University, 
      Princeton, NJ 08544, USA.}
      \quad \quad 
      Ricardo Masini\thanks{rmasini@ucdavis.edu, 
      UC Davis, 
      Davis, CA 95616, USA.}
      \quad \quad 
      Paul Schneider\thanks{paul.schneider@usi.ch, 
      USI Lugano and SFI, 
      Lugano, 6900, Switzerland.}}


\maketitle

\begin{abstract}
We develop a non-negative polynomial minimum-norm likelihood ratio (PLR) of two distributions of which only moments are known. The sample PLR converges to the unknown population PLR under mild conditions. The methodology allows for additional shape restrictions, as we illustrate with two empirical applications. The first develops a PLR for the unknown transition density of a jump-diffusion process, while the second extracts a positive density directly from option prices. In both cases, we show the importance of implementing the non-negativity restriction. 
\end{abstract}

\noindent%


\section{Introduction}\label{sec:intro}

In Hilbert spaces, orthogonality and minimum-norm problems are tightly related. As a consequence, orthogonal polynomials have a prominent role in minimum-norm approximation of unknown likelihood ratios and numeric integration. However, extant approaches based on orthogonality only do not preserve important properties of the approximated objects. In this paper, we develop projections of likelihood ratios onto positive polynomials, thus \emph{preserving positivity}, and if desired, additional structural constraints.

The literature considers approximations of likelihood ratios with orthogonal polynomials foremost due to the link between polynomials and the possibility of expressing moments of a distribution as expectations of polynomials \citep{aitsahalia02,filipovicmayerhoferschneider13,katokuriki13,rennerschmedders15}. While many other approaches exist in the literature to approximate likelihood ratios,\footnote{The literature on density approximations ranges from saddlepoint approximations \citep{aitsahaliayu05}, small time expansions \citep{yu07}, to simulation-based methods \citep{mijatovicschneider07,gieseckeschwenkler18} to cite a small subset of the econometrics and statistics literature. In machine learning,  starting from \citet{berlinetthomas-agnan04} a sequence of papers embeds distributions in  Reproducing Kernel Hilbert Spaces \citep{songetal09,gruenewaelderetal12,parkmuandet20,klebanovschustersullivan20}, as well as likelihood ratios \citep{schusteretal20}. While some of the aforementioned papers  preserve positivity and normalization of distributions, they do not provide the flexibility to include additional shape constraints.}  polynomials are a natural choice, as moments are known for many models such as the large affine class \citep{duffiefilipovicschachermayer03}. Moreover, under mild technical conditions, polynomials generate weighted $L^2$ spaces that arise naturally when working with expectations and sample averages. In our paper, we work with positive polynomials with the smallest modification to the extant expansions mentioned above.\footnote{More generally, we operate within the generic Hilbert space   problem  
\begin{align*}
\underset{g\in \mathcal H }{\text{ minimize }}  \norm{g} ,\text{ subject to } (g , f) \in C,  \, f \in \mathcal H
\end{align*}
where $C$ is a convex set, and  $\mathcal H$ is a Hilbert space along with its inner product $(\cdot,\cdot)$, encompassing function approximation, interpolation, and many other applications. In this paper, we specialize in likelihood ratios with probabilistic models in mind.} 

Our framework, rooted in \citet{grenander81} sieve estimation, is built upon several steps. First, we consider the conventional projection of a likelihood ratio on polynomials.
Through an optimization problem, we then obtain our polynomial minimum-norm likelihood ratio (PLR) as the sum of this element, and the minimum-norm polynomial that guarantees pointwise positivity. We show that this optimization problem, if feasible, has the PLR as its unique solution, and subsequently prove consistency 
of the PLR based on sample moments. Importantly, these coefficients of the PLR can be obtained rapidly as the solution of a conic optimization program, allowing also for additional constraints that can modify the shape of the PLR. 

The sieves literature nonparametrically approximates a function using a sequence of finite-dimensional spaces and corresponding basis functions, such as polynomials, splines, wavelets, and trigonometric functions,  that grow with the sample size to asymptotically eliminate the bias in the estimation process (see \cite{HandbookChenSieves}). In contrast, we opt to keep the degree of the PLR polynomial fixed, approximating the pseudo-true likelihood with (potential) bias. Our motivation stems from the fact that preserving positivity is more important when a finite sample size caps the maximum degree of the PLR than in a limiting large sample/large polynomial degree environment, where the estimated PLR would converge to its positive population counterpart.


A related literature beginning with \cite{hansenscheinkman09} emphasizes the importance of positive eigenfunctions in asset pricing problems\footnote{\cite{hansenscheinkman09} propose an operator approach to study the long-run risk-return trade-off properties of Stochastic Discount Factors (SDFs) in Markovian environments by solving a Perron-Frobenious eigenfunction problem. The SDF is decomposed into a permanent (martingale) and a transitory component (see \cite{alvarezjermann05}). \cite{ross15} uses the Perron-Frobenious approach to recover investors' beliefs, i.e., subjective probabilities under rational expectations.   \cite{borovickahansenscheinkman16} discuss the conditions under which the probability recovered with Perron-Frobenious coincides with the actual probabilities of investors under rational expectations. \cite{christensen17} proposes using Sieves to empirically identify this SDF decomposition, establishing the asymptotic theory for Perron-Frobenious eigenfunction and eigenvalue estimators. A non-exhaustive list of additional papers on the SDF decomposition and on recovery probabilities include \cite{bakshichabiyo12}, \cite{bakshichabiyogao15},  \cite{backuschernovzin14}, \cite{linetzkyqin16}, and \cite{linetzkyqin17}.}. Preserving positivity is fundamental in such problems, for instance, in \cite{christensen17}, who provide Hermite polynomial approximations for estimating a nonparametric decomposition of the Stochastic Discount Factor in transitory and permanent components. In this context, our methodology could be used to guarantee a non-negative approximation to the positive eigenfunctions that define the SDF decomposition in finite samples where a projection on polynomials is not guaranteed to preserve positivity.
     
We illustrate the usefulness of PLR with two simple applications. In the first, we expand the likelihood ratio of the transition density of a continuous-time jump-diffusion process with respect to a Gamma density.\footnote{When expanding a transition density it is necessary to choose an appropriate auxiliary density function on the same support. For positive support, the Gamma distribution is suited (see \citet{filipovicmayerhoferschneider13}).} While the conventional orthogonal polynomial approach produces an approximation that is negative close to zero, our PLR is non-negative everywhere and thus can be readily used in any context that crucially relies on 
the structural properties, normalization, and positivity, of probability measures, such as accept-reject sampling or Bayesian modeling. 

In our second application, we subsequently investigate how to extract probability densities directly from option prices. This is particularly interesting in the context of pricing exotic and/or over-the-counter derivatives based on the probability density extracted from vanilla observed option prices. We illustrate this procedure in a controlled environment with a bivariate continuous-time jump-diffusion process from \citet{duffiepansingleton00}. Our experiments reveal how the 
positivity of the PLR crucially translates into meaningful estimates of positive probability densities, rather than projections of densities that do not share basic defining properties of probability measures.

The paper is organized as follows. In Section \ref{sec:setup} we develop the necessary notation. The traditional approach to minimum-norm expansions appears in Section \ref{sec:rkhspzn} for reference. In Section \ref{sec:genprob} we propose modifications to the standard program to ensure positivity. Section \ref{sec:properties} shows the uniqueness, consistency, and rate of convergence of the polynomial coefficients defining the PLR. Applications are presented in Sections \ref{sec:densprox} (density approximations), and \ref{sec:optionpricing} (option pricing).
Section \ref{sec:conclusion} concludes. In the  \hyperref[sec:pospol]{Appendix}, we review sum-of-squares (s.o.s.) polynomials and include the proofs for the theoretical results. 

\section{Constrained polynomial likelihood}\label{sec:plr}
\subsection{Set-up and notation}\label{sec:setup}
Denoting by $\PP$ a probability measure on $D\subseteq \reals ^{d}$, let  $L_{\PP}^{2}$ be the equivalence class of functions $f:D\mapsto \mathbb \reals$  such that $\int _{D}f^2(t)d\PP(t)<\infty$ with inner product
\begin{equation}\label{eq:innerprod}
\begin{split}
 (f,g)=\int _D {f(t)g(t)d\PP( t)} , \, f,g\in L^{2}_\PP, \,  t:=[t_1,\ldots, t_d],
\end{split}
 \end{equation}
 where we allow in subsequent sections to extend the inner product to act element-wise on matrices. 

Denote by $\reals[t]$ the ring of square-integrable polynomials on $\reals ^{d}$, and by $\reals [t]_n$  the subset of polynomials  $\xi \in \reals [t]$ with  $\deg (\xi)\leq n$. Denote further by $P_{\PP,n}$ the set of polynomials $\reals [ t]_n$  endowed with inner product \eqref{eq:innerprod},  a finite-dimensional Hilbert space.

From  the standard canonical monomial basis of $P_{\PP,n}$ we denote
\begin{equation}\label{eq:monombasis}
\bm \tau _n( t):=[1, t_1,\ldots, t_d, t_1^2,t_1 t_2,\ldots, t_d^{2},\ldots,t_1^{n},t_1^{n-1}t_2,\ldots, t_d^{n}]^{\T},
\end{equation}
as well as  multi-index powers $t^{\bm \beta}:=t_1^{\beta _{1}}\cdots t_d^{\beta _{d}}$ for $\bm \beta \in \mathbb N^{d}_0$, where  the length of the multi-index is $\abs{\bm \beta}=\beta_1+\cdots+\beta _d$ . 

There are $\binom{n+d}{d}$ elements in the monomial basis, and we denote by $\alpha_0,\ldots, \alpha _{N}$, with  $N=\binom{n+d}{d}-1$, the  multi-indices corresponding to their order of appearance in   \eqref{eq:monombasis}. For example, the second element  is $t^{\bm \alpha_{1}}=t_1^{\alpha _{11}}\cdots t_d^{\alpha _{1d}}=t_1^1t_2^0\cdots t_d^0=t_1$.  Any  polynomial $\xi_n\in P_{\PP,n}$ can  be written as $\xi_n=\bm x^{\T}  \bm \tau_n$, where $\bm x\in \reals ^{N+1}$ is a coefficient vector. 

With projections on positive polynomials in mind, we parameterize a particularly tractable cone $M_n(D)\subset P_{\PP,n}$ of positive polynomials as follows.
\begin{definition}[Positive s.o.s. polynomials]\label{def:MD}
Let $S_+^q, q\in \mathbb N$ denote the set of symmetric positive semidefinite (p.s.d.) matrices of dimension $q$. We denote by $M_n(D)$ the set of positive polynomials $\xi_n(t)=\bm x^{\T}\bm \tau _n(t)\geq 0, \forall t \in D$ for matrices $\bm V\in S_+^{q_V}$ and $\bm W\in S_+^{q_W}$ on a case-by-case basis.
\begin{enumerate}
 \item Case $D= \reals^{d}$: $\xi_n(t)=\bm \tau_{n/2}^{\T}(t) \, \bm V \, \bm \tau_{n/2}(t), $
 \item Case $D=\reals _+$: $\xi_n(t)=\bm \tau_{n/2}^{\T}(t) \,\bm V \, \bm \tau_{n/2}(t)+t \brac{\bm \tau_{n/2-1}^{\T}(t) \, \bm W \, \bm \tau_{n/2-1}(t)},$,
 \item Case $D=[a,b]$: $\xi_n(t)=\bm \tau_{n/2}^{\T}(t) \, \bm V \, \bm \tau_{n/2}(t)+(b-t)(t-a) \brac{\bm \tau_{n/2-1}^{\T}(t) \, \bm W \, \bm \tau_{n/2-1}(t)}.$
\end{enumerate}
 \end{definition}
 By collecting coefficients in the relation above, we can identify a linear relation $\bm x=T(\bm V, \bm W)$ between the coefficient vector $\bm x$ of a positive polynomial in $M_n(D)$  to elements of the symmetric positive semidefinite (p.s.d.) matrices $\bm V$ and $\bm W$, whose dimensions are $q_V=\binom{n/2+d}{d}$ and $q_W=0$ in the first case, and $q_V=\binom{n/2+1}{1}$ and $q_W=\binom{n/2}{1}$ in the latter two cases. As outlined in Appendix \ref{sec:pospol}, Definition \ref{def:MD}  characterizes the set of positive polynomials in the univariate case, but it is only sufficient in the multivariate case, for which positive polynomials exist that are not sums of squares. As the benefit of the above s.o.s. parameterization,  positive polynomials can be modellied through semidefinite programming. 
 Without further modifications, the linear transformation $T$ is not a bijection in general, however, as the following example shows.

 \begin{example}[Selection matrices]\label{ex:selection}
  Suppose that $n=4$, and $D=\reals$. Then positivity requires
  \begin{equation*}
   \bm \tau_4^{\T}(t)\bm x = \bm \tau_{2}^{\T}(t) \, \bm V \, \bm \tau_{2}(t)=\tr \brac{\bm V \bm \tau_{2}(t)\bm \tau_{2}^{\T}(t) }=\tr \brac{\begin{bmatrix}
V_{11} & V_{21} & V_{31} \\
V_{21} & V_{22} & V_{32} \\
V_{31} & V_{32} & V_{33}
\end{bmatrix}
 \begin{bmatrix}
               1 & t & t^2 \\                                                                                                                                      
               t & t^2 & t^3 \\
               t^2 & t^3 & t^4
               \end{bmatrix}
}.
  \end{equation*}
Collecting coefficients, we can thus express the selection matrices as 
\begin{equation*}
 \bm V_0=\begin{bmatrix}
          1 & 0 & 0 \\
          0 & 0 & 0 \\
          0 & 0 & 0
         \end{bmatrix}, 
 \bm V_1=\begin{bmatrix}
          0 & 1 & 0 \\
          1 & 0 & 0 \\
          0 & 0 & 0
         \end{bmatrix}, 
          \bm V_2=\begin{bmatrix}
          0 & 0 & 1 \\
          0 & 1 & 0 \\
          1 & 0 & 0
         \end{bmatrix}, 
         \bm V_3=\begin{bmatrix}
          0 & 0 & 0 \\
          0 & 0 & 1 \\
          0 & 1 & 0
         \end{bmatrix}, 
          \bm V_4=\begin{bmatrix}
          0 & 0 & 0 \\
          0 & 0 & 0 \\
          0 & 0 & 1
         \end{bmatrix},
\end{equation*}
such that $[  x_0, x_1, x_2, x_3, x_4]=
  T^{\T}(\bm V)=[\tr \brac{\bm V\bm V_i}]_{i=0}^4=[ V_{11}, 2 V_{21}, V_{22}+2V_{31}, 2V_{32},  V_{33}].$
 \end{example}

Before presenting the full type of problems we are intending to solve on $M_n(D)$, we briefly review projections of Radon-Nikodym derivatives on polynomials through equality-constrained minimum-norm problems. While this problem is standard, and its solution is well-known, it helps inform the optimization program for the shape-constrained version in the subsequent section. 

\subsection{Minimum-norm likelihood ratio projection}\label{sec:rkhspzn}
For projections in  $L^2_\Pmeas$, with or without shape constraints, the   Gram matrix
\begin{equation}\label{eq:momentmatrix}
 \bm H_n:=(\bm \tau _n, \bm \tau _n ^{\T})=\begin{bmatrix}
      \mu _{0,0}^\PP & \mu _{0,1}^\PP & \cdots & \mu _{0,N}^\PP \\
      \mu _{1,0}^\PP & \mu _{1,1}^\PP & \cdots & \mu _{1,N}^\PP \\
      \vdots & \vdots & \ddots & \vdots \\
      \mu _{N,0}^\PP & \mu _{N,1}^\PP & \cdots & \mu _{N,N}^\PP
     \end{bmatrix},\text{ with } \mu_{i,j}^\PP:=\int _D t^{\alpha _i+\alpha _j}d\PP(t)
\end{equation}
 features prominently.  Our first assumption regarding this symmetric and positive semidefinite matrix greatly simplifies further analysis. It is  mild, easily checked, and
ensures in particular that the inverse of the Gram matrix $\bm H^{-1}$ is readily available, greatly facilitating the analysis below.
 \begin{assumption}[Rank of $\bm H_n$]\label{ass:card}Gram matrix $\bm H_n$'s smallest eigenvalue $\lambda _{\min}(\bm H_n)>0$.
\end{assumption}
We are particularly interested in projections on polynomials of the Radon-Nikodym derivative of a probability measure $\QQ$ on $D$ that is absolutely continuous with respect to $\PP$.  Finding the minimum-norm such projection on polynomials, where $\QQ$ is only known from its moment vector $\bm \mu _\QQ:=[\mu_0^\QQ,\ldots, \mu_N^\QQ]^{\T}$ with $\mu_{i}^\QQ:=\int _D t^{\alpha _i}d\QQ( t)$,  solves the  equivalent problems \citep{luenberger97}

\noindent\begin{minipage}{.45\textwidth}
\begin{equation*}
 \underset{\xi_n \in P_{\PP,n}}{\text{minimize }}  \norm{\frac{d\QQ}{d\PP}-\xi _n}^2 
\end{equation*}
\end{minipage}
$\Leftrightarrow$
\begin{minipage}{.45\textwidth}
\begin{equation}\label{eq:thepznprogram}
\begin{split}
 \underset{\bm x \in \reals ^N}{\text{minimize }} & \norm{\bm H_n^{1/2} \, \bm x}_2^2, \, \text{subject to } \\
  \bm H_n \bm x&=\bm \mu _\QQ,  \\
  \end{split}
\end{equation}
\end{minipage}
where $\bm H _n^{1/2}$ denotes the Cholesky factor of $\bm H_n$.  The well-known  solution to problem \eqref{eq:thepznprogram} is
\begin{equation}\label{eq:fms}
 \xi_n^{\star}(t):=\bm \mu _\QQ^{\T}\bm H_n^{-1}\bm \tau _n( t).
\end{equation}

On $P_{\PP,n}$, the projection of a Radon-Nikodym derivative on polynomials evidently  corresponds to moment-matching. However, there is no mechanism that ensures to maintain the structural properties of Radon-Nikodym derivatives  -- normalization and positivity -- in the projection. 

To project on \emph{positive} polynomials, with  $M_n(D)$ not being a linear space, using it as a hypothesis space requires additional considerations.  In the next section, we, therefore, introduce an optimization program akin to \eqref{eq:thepznprogram}, where membership to $M_n(D)$ is implemented as a conic constraint. The formulation as an optimization program additionally allows to generalize the equality constraints in Section \ref{sec:rkhspzn}  to inequalities, and to add further restrictions, as demanded by each application.

\subsection{Constrained minimum-norm likelihood ratio projection}\label{sec:genprob}
In this section, we introduce a minimization program that accommodates the constraints in  \eqref{eq:thepznprogram} as a special case and includes additionally the conic constraint $\xi_n\in M_n(D)$ from Definition \ref{def:MD}. 
\begin{equation}\label{eq:theprogram}
\begin{split}
 \underset{  \xi_n \in M_n(D)}{\text{minimize}} & \norm{\xi_n}^2,  \text{ subject to } \\
  (\xi_n , f_i)&=c_i, \, i=1,\ldots,m, \text{ and }   (\xi_n , g_j)\leq d_j, \, j=1,\ldots,l,
  \end{split}
\end{equation}
where $f_1,\ldots, f_m, g_1,\ldots, g_l\in L^{2}_\PP$ are linearly independent  functions generating the subspace $K:=[f_1,\ldots, f_m, g_1,\ldots, g_l]$. While in general  $K\varsubsetneq P_{\PP,n}$,    
we nevertheless routinely set $f_1=1$, so that with   $(\xi _n,f_1)=1$, $\xi _n$ represents a normalized and positive likelihood ratio with respect to a probability distribution $\PP$. Together with the additional constraints it therefore is a constrained polynomial likelihood ratio (PLR). Note that the constraint $\xi_n\in M_n(D)$ ensures pointwise positivity. This is important for evaluating the PLR with any argument while maintaining positivity, as for instance in out-of-sample exercises or simulations. 
    
It is easy to construct an example for which program \eqref{eq:theprogram} is infeasible. The constraint $(\xi_n,1)=-1$ can not be satisfied by any non-negative polynomial $\xi_n$,  with $\PP$ being a probability measure. Apart from such cases, the feasibility of \eqref{eq:theprogram} becomes more likely, the higher the order of $\xi_n$, increasing the number of coefficients. In light of this observation, denoting by $\tilde K =\{\bm x \in \reals ^{N+1} : \{(\xi_n , f_i)=c_i\}_{\{i=1,\ldots,m\}}, \text{ and } \{(\xi_n , g_j)\leq d_j\}_ {\{j=1,\ldots,l\}}\}$ the feasible set of program \eqref{eq:theprogram}, we state our second assertion.
 \begin{assumption}[Feasibility]\label{ass:feasible}  There exists $n\in \mathbb N$  such that  $\tilde K$  has non-empty interior. 
\end{assumption}

In the next section, we discuss a generic solution algorithm to the optimization problem \eqref{eq:theprogram}, and show the consistency of the solution in the case when the inner product \eqref{eq:innerprod} is estimated from sample averages.
    
\subsection{Properties of solutions on $P_{\PP,n}$ and $M_n(D)$}\label{sec:properties}
The present section discusses the properties of the solution to program \eqref{eq:theprogram}.  

\subsubsection{Uniqueness of the primal optimization problem}

As a first step, we argue the uniqueness of the solution to program \eqref{eq:theprogram}, if the constraints are feasible, as a standard result in finite-dimensional convex optimization. To develop a solution, we first exploit that we work in a finite-dimensional Hilbert space of polynomials, which allows us to express the functional inequalities as matrix equations. The (in)equalities in Eq. \eqref{eq:theprogram} are linear in  $\xi_n$, and defining $\bm c :=[c_1,\ldots,c_m] ^{\T}, \, \bm d:=[d_1,\ldots,d_l]^{\T},\, \bm f:=[f_1,\ldots, f_m]^{\T}$, and $\bm g:=[g_1,\ldots, g_l]^{\T}$, we can express them as 
\begin{equation}\label{eq:constraintscoordinate}
 (\xi_n,\bm f)=\bm F_n \bm x=\bm c, \text{ and }
 (\xi_n,\bm g)=\bm G_n \bm x\preceq \bm d.
\end{equation}
where $\preceq$ represents generalized, conic, inequality, and the inner product is naturally applied element-wise, such that $\bm F_n\in \reals ^{m\times N+1}$ and $\bm G_n\in \reals ^{l\times N+1}$.

\begin{proposition}\label{prop:unique}
 If Assumptions \ref{ass:card} and \ref{ass:feasible} hold, problem \eqref{eq:theprogram} has a unique solution from the  primal problem
 \begin{equation}\label{eq:mixedconicprimal}
\begin{split}
\underset{\bm x\in \reals ^{N+1}}{\text{minimize  }}  & \quad \frac{1}{2}\norm{\bm H_n^{1/2}\bm x}_2^2,   \\
 \text{subject to } & 
   \bm F_n \bm x=\bm c, \, 
   \bm G_n \bm x\preceq  \bm d , \\
    \bm x &=T(\bm V, \bm W);  
   \bm V\in S_+^{m_V}, \bm W \in S_+^{m_W}, 
   \end{split}
\end{equation}
where the linear map $T:S_+^{m_V}\times S_+^{m_W}\mapsto \reals ^{N+1}$,
\begin{equation}\label{eq:linop}
 T(\bm V, \bm W):=[
                \tr (\bm V_{0} \bm V)+\tr (\bm W_{0} \bm W) , 
                \cdots ,
                \tr (\bm V_{N} \bm V)+\tr (\bm W_{N} \bm W)
          ] ^{\T},
\end{equation}
for fixed symmetric selection matrices $\bm V_{0},\ldots, \bm V _{N}$ and $\bm W_{0},\ldots, \bm W _{N}$ of the same dimension as $\bm V$ and $\bm W$, respectively. 
The corresponding dual reads
 \begin{equation}\label{eq:mixedconicdual}
\begin{split}
 \underset{\bm \eta \in \reals ^m, \bm \varepsilon \in \reals _+^l, \bm \nu \in \reals ^{N+1}}{\text{maximize }} & -\frac{1}{2}\|\bm H^{-1/2}_n (\bm F_n ^{\T}\bm \eta + \bm \nu -\bm G _n^{\T}\bm \varepsilon) \|_2^2-\bm \varepsilon ^{\T}\bm d+\bm \eta ^{\T} \bm c\\
 \text{subject to } &\sum _{i=0}^N\nu_i\bm V_i  \succeq  0,\sum _{i=0}^N\nu_i\bm W_i  \succeq  0.
   \end{split}
\end{equation}
where   the primal and dual solutions are linked  as
\begin{equation}\label{eq:myfoc}
 \bm x_0 = \bm H_n ^{-1}(\bm F_n ^{\T}\bm \eta _0 + \bm \nu _0 -\bm G _n^{\T}\bm \varepsilon _0).
\end{equation}
\end{proposition}
  It is noteworthy that the dual formulation is strictly convex in the sum $\bm F_n ^{\T}\bm \eta + \bm \nu -\bm G _n^{\T}\bm \varepsilon$, rather than in the Lagrange multipliers $\bm \eta, \bm \nu, \bm \varepsilon$ corresponding to the equality, positivity, and inequality constraints, separately. This observation becomes important for deriving results about the behavior of optimization problems \eqref{eq:mixedconicprimal} and \eqref{eq:mixedconicdual} when the  population matrices $\bm F_n,\bm G_n, \bm H_n$, and possibly also $\bm c$ and $\bm d$ are estimated from data. Conveniently,  the ingredients of the program \eqref{eq:theprogram} depend entirely and exclusively on  expectations induced by $\PP$. This feature is useful in applications where the weight function $\PP$ in \eqref{eq:innerprod} is unknown, but, for instance,  its moments can be estimated.

  Furthermore, note that if $K$ is generated by polynomials,  we have a unique orthogonal decomposition of  solution $\xi_n(t)=\bm x_0 ^\T \bm \tau_n(t)$ of program \eqref{eq:mixedconicprimal}  into
  \begin{equation}\label{eq:decomposition}
  \xi_n=\xi_n^{\star}\oplus \xi_n^{\circ},
  \end{equation}
where $\xi_n^{\star}$ is the minimum-norm polynomial on $K$, and $\xi_n^{\circ}$ is the minimum-norm polynomial lifting the solution into $M_n(D)$.
\subsubsection{Consistency}\label{sec:consistency}
In this section, we are concerned with the behavior of the optimal coefficients when the inner product can only be evaluated from the sample distribution 
induced by identically distributed draws $x_1,\ldots, x_k$ from $\PP$. 
We derive two types of asymptotic results, the first being a consistency result on the coefficients $\bm x$, and the second being a non-standard result accounting for the non-uniqueness of the arguments that minimize the dual \eqref{eq:mixedconicdual}.

We first define 
\begin{equation}\label{eq:LLN}
 \langle f,g \rangle _k:= \frac{1}{k}\sum _{i=1}^kf(x_i)g(x_i) \text{ for }f,g\in L^{2}_\PP,
\end{equation}

With the empirical inner product, we can define  $\hat {\bm H}_{n,k}:=\left \langle \bm \tau _{n}, \bm \tau _n^{\T}\right \rangle_k$, and analogously
 $\hat {\bm F}_{n,k}:=\langle \xi_n, \bm f \rangle _k$ and $\hat {\bm G}_{n,k}:=\langle \xi_n, \bm g \rangle _k$. Denoting by $\tilde K_k$ the feasible set induced by $\hat {\bm F}_{n,k},\hat {\bm G}_{n,k}$, and $\hat {\bm H}_{n,k}$, we are in the position to introduce the sample versions of problem \eqref{eq:theprogram}, of its corresponding mixed-conic primal and dual matrix versions \eqref{eq:mixedconicprimal} and \eqref{eq:mixedconicdual}, and to establish consistency of the estimators in these problems.
Let $\bm x _0$ represent the coefficients of the polynomial $\xi_n$, that are the solution of problem \eqref{eq:theprogram} with restrictions based on the population moments $\{(\xi_n,f_i),(\xi_n,g_j)\}$, and $\hat{\bm x}_k$ represent the corresponding coefficients from the sample moments $\{\langle \xi_n,f_i\rangle _k,\langle \xi_n,g_j\rangle_k \}, i=1,\ldots,m; j=1,\ldots,l$, with analogous notation for the dual problem.

To provide the two types of consistency results,  we introduce the population and sample objective functions
\begin{equation*}
 \RR_k(\bm x):=-\bm x^{\T}\hat{\bm H}_{n,k}\bm x \text{ and } \RR_0(\bm x):=-\bm x^{\T}\bm H_n \bm x \text{ of the primal problem.}
\end{equation*} 
 
 To establish consistency of the estimators, as well as statements about their behavior with increasing $k$, we make a (mild) assumption concerning the  matrices depending on estimated moments.
 \begin{assumption}[Rate of Convergence]\label{ass:rate}
 Let $\widehat{\bm H}_{n,k}$, $\widehat{\bm F}_{n,k}$, $\widehat{\bm G}_{n,k}$, $\widehat{\bm c}_k$ and $\widehat{\bm d}_k$ be estimators of $\bm H_n$, $\bm F_n$, $\bm G_n$, $\bm c$ and $\bm d$, respectively. Set $\bm s^{\T} :=[(\vc \ \bm H_n)^{\T},(\vc\ \bm F_n)^{\T},(\vc\ \bm G_n)^{\T}, \bm c ^{\T}, \bm d ^{\T}]$ and define $\widehat{\bm s}$ in the same manner. Suppose that $\|\widehat{\bm s} - \bm s\|_2=O_\PP(r_k)$  for some non-negative sequence $r_k\to 0$ as $k\to\infty$.
\end{assumption}

In most cases, the estimators $\widehat{\bm H}_{n,k}$, $\widehat{\bm F}_{n,k}$, $\widehat{\bm G}_{n,k}$, $\widehat{\bm c}_k$ and $\widehat{\bm d}_k$ will be sample averages. Since we already assume square integrability, the Weak Law of Large numbers will ensure Assumption \ref{ass:rate} with $r_k=1/\sqrt{k}$. Uncorrelatedness of the realizations suffices, but several other forms of weak dependencies could be accommodated such as mixingales and mixing sequences provided that the serial correlation can be controlled. Note that Assumption \ref{ass:rate} together with Assumption \ref{ass:card} also ensures that the symmetric p.s.d. matrix $\hat {\bm H}_{n,k}$ will be positive definite with high probability for $k$ large enough.

\begin{proposition}[Consistency]\label{pro:consistency}
Under Assumptions \ref{ass:card}, \ref{ass:feasible}, and \ref{ass:rate}, $\hat{\bm x}_k \overset{\PP}{\to} \bm x_0$ as $k \to \infty$ for  fixed   $n$.
\end{proposition}

Proposition \ref{pro:consistency} is based on mild assumptions and describes the asymptotic behavior of the estimator of the optimal coefficients $\bm x _0$. However, it does not account for the p.s.d. matrices $\bm \Lambda _V$ and $\bm \Lambda _W$, which appear as Lagrange multipliers of the s.o.s. cones (cf. the Lagrangian in \eqref{eq:mylagrangian}). 

To describe the joint behavior of the estimated primal and dual parameters,  let $\bm\theta^{\T}:=\big[\bm x^{\T},(\vecme \bm V)^{\T},(\vecme \bm W)^{\T},\bm\eta^{\T},\bm\varepsilon^{\T},\bm\nu^{\T},(\vecme\bm \Lambda_V)^{\T},(\vecme \bm \Lambda_W)^{\T}\big]\in \Theta$ where $\Theta$ is the Cartesian product of the domain of each component of $\bm\theta$ defined above, and $\Psi:\Theta\to\mathbb{R}^e$ is given by
\[
\Psi(\bm\theta) :=\begin{bmatrix}
    \bm H _n \bm x  - \bm F^{\T}_n\bm\eta + \bm G^{\T}_n\bm\varepsilon -\bm\nu\\
    \vc(\bm \Lambda_V -\sum_{i=0}^n\nu_i \bm V_i)\\
     \vc(\bm \Lambda_W -\sum_{i=0}^n\nu_i \bm W_i)\\
     \bm F_n \bm x  -\bm c\\
     \bm\varepsilon\odot(\bm G _n\bm x  -\bm d)\\
     \bm x - T(\bm V,\bm W)\\
     \tr(\bm \Lambda_V\bm V)\\
     \tr(\bm \Lambda_W\bm W)\\
\end{bmatrix},
\]
where $\odot$ denotes the Hadamard product. The solutions to problems \eqref{eq:mixedconicprimal}--\eqref{eq:mixedconicdual} can be characterized by $\Theta_0:=\{\bm\theta \in \Theta: \Psi(\bm\theta)=\bm 0\}$ due to strong duality. Note that $\Psi$ is closely related to the KKT conditions, in the sense that the first three equations of $\Psi(\bm\theta)=\bm 0$ are stationary conditions, and the remaining are non-slackness-type, and primal and dual feasibility conditions.

Since $\Theta_0$ is a level set of a continuous function, it is closed. It is also convex as the solution set of a convex optimization problem, but it is not necessarily bounded. 
However, for optimal $\bm x _0\in\mathbb{R}^{N+1}$, we can regularize such that all other entries are bounded by a (large) constant $C>0$. Take, for instance, Example \ref{ex:selection}, where $x_2=V_{22}+2V_{31}$. We can impose a large upper bound $C$ for the magnitudes of $V_{22}$ and $V_{31}$. Following this principle,  we consider only solutions in the compact $\Theta_*:= \Theta_0\cap \bar{B}_C(\bm {0})$, where $\bar{B}_C(\bm {0})$ denotes the closed ball in $\Theta$ with radius $C$. 

Let $\widehat{\Psi}$ be the sample analogue of  the function $\Psi$, i.e., $\widehat{\Psi}$ is the function $\Psi$ with $\bm H _n$, $\bm F _n$, $\bm G_n$, $\bm c$ and $\bm d$ replaced by estimators $\widehat{\bm H}_{n,k}$, $\widehat{\bm F}_{n,k}$, $\widehat{\bm G}_{n,k}$, $\widehat{\bm c}_{k}$ and $\widehat{\bm d}_{k}$, respectively, and define $\widehat{\Theta}_k:=\{\bm\theta \in  \bar{B}_C(\bm {0}):\|\widehat{\Psi}(\bm\theta)\|\leq \delta_k\}$ for some non-negative sequence such that $\delta_k\to 0$ as $k\to\infty$, keeping in mind Assumption \ref{ass:rate} for the estimators.

Due to the partial identification of all coefficients (except for $\bm x$ that is unique), both the estimator $\widehat{\Theta}_k$ and the parameter of interest $\Theta_*$ are sets (subsets of $\Theta$). Therefore, we use the Hausdorff distance defined by a pair of sets of a metric space $A,B\subseteq (M,d)$ as
\[
d_H(A,B):=\max\left\{\sup_{a\in A}\inf_{b\in B} d(a,b), \sup_{b\in B}\inf_{a\in A} d(a,b)\right\}.
\]

\begin{proposition}[Set rate of convergence]\label{pro:rate} Under Assumptions \ref{ass:feasible} and \ref{ass:rate}, for every $\epsilon>0$, there are constants $C_\epsilon, M_\epsilon>0$ such that for every $k\in \mathbb{N}$,
\[\PP(d_H(\widehat{\Theta}_k,\Theta_*)\leq M_{\epsilon}r_k)\geq 1-\epsilon,
\]
provided that $\delta_k =  C r_k$ for $C
\geq C_\epsilon$ where $d_H$ is the Hausdorff distance.
\end{proposition}

Proposition \ref{pro:rate} above complements Proposition \ref{pro:consistency} by establishing the rate of convergence of $\widehat{\Theta}_k$ to $\Theta_*$ as $k\to\infty$ in the Hausdorff distance, which coincides with the usual Euclidean distance when both $\widehat{\Theta}_k$ and $\Theta_*$ are singletons. In general, we expect $r_k = 1/\sqrt{k}$, and in that case, $\widehat{\Theta}_k$ is $\sqrt{k}$-consistent, which implies that $\widehat{\bm x}_k - \bm x_0= O_\PP(1/\sqrt{k})$. From the proof of Proposition \ref{pro:rate}, the following related result follows immediately. 

\begin{corollary}[Asymptotic coverage]\label{cor:coverage} Under Assumptions \ref{ass:feasible} and \ref{ass:rate} we have $\PP(\Theta_*\subseteq \widehat{\Theta}_k)\to 1$ provided that $\delta_k =O(r_k \log k)$ as $k\to\infty$. 
\end{corollary}

Corollary \ref{cor:coverage} provides a useful guarantee for the estimated coefficients of the model. It ensures that the $(r_k\log k)$-level set of $\widehat{\Psi}$ contains the true set $\Theta_*$ with high probability as $k$ increases. Once again, we expect $r_k = 1/\sqrt{k}$, consequently $\widehat{\Theta}_k$ is a (nearly parametric rate) confidence set for $\Theta_*$ with confidence level approaching one as the sample size increases. Note that the $\log$ term is arbitrary in the sense that any strictly increasing function could replace it without changing the result.

\section{Applications}

In this section, we present two applications of the PLR developed in this paper. The first, in Section \ref{sec:densprox}, approximates the unknown transition density of a jump-diffusion process  whose moments are known. The second extracts densities from option prices in Section \ref{sec:optionpricing}. For both applications,  both $\PP$, as well as $\QQ$ possess Lebesgue densities $p$ and $q$, respectively. 

\subsection{Density expansions}\label{sec:densprox}
As a natural application and as a continuation of Section \ref{sec:rkhspzn}, we investigate here approximations of an unknown distribution $\QQ$ given moments $\mu _\QQ :=\int _D \tau _md\QQ$, where we allow for $m\leq n$ (rather than $m=n$ as in Section \ref{sec:rkhspzn}). 
To indicate the dependence of $\xi_n$ on the number of monomials spanning $K$,  we denote the optimal $n-$order polynomial likelihood expansion with moment constraints up to order $m$ by $\xi ^{(m)}_{n(m)}$ or more concisely $\xi ^{(m)}_{n}$, where $n$ is chosen sufficiently high according to Assumption \ref{ass:feasible}. 
We begin by showing that $\xi  ^{(m)}_n$ converges in $L^{2}_\PP$ to its unknown counterpart $\frac{d\QQ}{d\PP}$. For this asymptotic analysis formulated with population moments, we need a  technical assumption.\footnote{It is important to keep in mind that for fixed $n$ the framework from Sections \ref{sec:rkhspzn} and \ref{sec:genprob} is applicable also to distributions with finite moments, but no moment-generating function, like the log-normal distribution. However, for the asymptotic analysis on the polynomial dimension $n$ as shown in this section, only $\PP$ random variables with a moment-generating function are permissible.} 
\begin{assumption}[Polynomial basis]\label{ass:basis}
 The ring of polynomials $\reals[t]$ is a basis of $L^{2}_\PP$.
\end{assumption}
Assumption \ref{ass:basis} is justified in particular for compact state spaces, as well as unbounded state spaces with the tails of $p$ decaying sufficiently quickly \citep{filipovicmayerhoferschneider13}. 

\begin{theorem}\label{th:L2}
If   $\frac{d\QQ}{d\PP}\in L^{2}_{\PP}$, and Assumptions \ref{ass:card}, \ref{ass:feasible}, and  \ref{ass:basis}  hold, the non-negative expansion $\xi^{(m)}_n$ converges in $L^2_{\PP}$,
\begin{equation}
 \lim _{m\to \infty}\norm{\frac{d\QQ}{d\PP} -\xi ^{(m)}_{n(m)}}=0.
\end{equation} 
\end{theorem}

To put our PLR to work, we now confront our density approximation approach with the one proposed in \citet[][FMS]{filipovicmayerhoferschneider13}. For this purpose, we consider the basic affine jump diffusion (BAJD) solving the stochastic differential equation   
\begin{equation}\label{eq:ajd}
dY_t=(\kappa \theta - \kappa  Y_t)\,dt+\sigma \sqrt{Y_t}dW_t +dL_t.
\end{equation}
where $W_t$ is a Brownian motion, the intensity of the compound Poisson process $L$ is $\lambda \geq 0$, and
the expected jump size of the exponentially distributed jumps is
$\nu \geq 0$. The transition density of  $Y_{\Delta}\mid Y_0, \, \Delta >0$ is not known in closed form, but its existence is assured if $2 \kappa \theta > \sigma ^2$ \citep[Theorem 2]{filipovicmayerhoferschneider13} on its domain $D=\reals _+$. Note that since the BAJD is a polynomial process, its conditional moments $\mu_i:=\mathbb E\brak{Y_{\Delta}^i\mid Y_0}$ are known in closed form for $\Delta>0$, however, even though the transition density is not. This process, as well as some of its variations, have been adopted to model the dynamics of stock index prices (\cite{bates00}), and to represent the intensity of the first jump to default when pricing CDS options and other related credit derivatives (\cite{brigomercurio01}). 

In the following, we develop a likelihood ratio tilting a Gamma distribution $\Gamma(1+\tilde{p},1)$ with density $p(x;\tilde{p})=\frac{e^{-x}x^{\tilde{p}}}{\Gamma\brak{1+\tilde{p}}}$, 
where $\tilde{p}=\mu_1^2/(\mu_2-\mu _1^2)-1$ and $\Gamma$ denotes the Gamma function. This auxiliary density pertains to the scaled random variable $\bar Y_\Delta:=Y_\Delta \cdot c$, with $c:=\frac{ \mu_1}{\mu_2-\mu_1^2}$, so that the object of interest $q^{(m)}(y)=\xi _n^{(m)}(\bar y) p(\bar y)\cdot c$ accounting for this change of variables, where the notation $\xi _n^{(m)}(\bar y)$ emphasizes the dependence of the PLR on the scaled random variable. 

We match moments $\mu_0,\ldots, \mu _5$, where the choice of $m=5$ is arbitrary, leaving sufficiently many free coefficients ($n=8$) to obtain a PLR. The corresponding program reads
\begin{align*}
 \underset{\xi_8\in M_8(\reals _+)}{\text{minimize }} & \norm{\xi_8},  \text{subject to }  
& (\xi_8,t^i)=\mu_i ,\, i=0,\ldots, 5,
\end{align*}
and we denote by $\xi _8^{(5)}$ the solution to the program above.
We confront our PLR approximation with the one proposed in FMS using the same weight function $p$, and matching the same moments $\mu _0,\ldots, \mu _5$. The corresponding program reads
\begin{align*}
 \underset{\eta_5\in \reals[t]_5}{\text{minimize }} & \norm{\eta_5}, \text{subject to } 
 & (\eta_5,t^i)=\mu_i ,\, i=0,\ldots, 5.
\end{align*}
and  we term its solution $\eta ^{(5)}$.  FMS solves the program using the projection theorem via orthogonal polynomials. Ours and their solutions are easily related, as $\xi ^{\star(m)}_{n}=\eta ^{(m)}$ for every non-negative integer $m$ by construction, so that in general $ \xi ^{(m)}_{n}=\eta ^{(m)}+\xi ^{\circ(m)}_{n}$ (cf. decomposition \eqref{eq:decomposition}).
From elementary arguments, the PLR is farther away in $L^{2}_\PP$ norm from the true likelihood ratio than the FMS one for every $m$. This is the price of non-negativity, which we will see in our next illustration might be worth paying, depending on the application.

We perform the comparison with the parameters $\kappa \theta = 0.05, \kappa = 1, \sigma = 0.2, \lambda=1, \nu=0.05, y_0=0.05$, roughly describing the dynamics of a stochastic equity volatility process,  with the true transition density obtained numerically from Fourier inversion using the exponentially-affine characteristic function of the BAJD. 

\begin{figure}
\centering
\begin{subfigure}[c]{0.46\textwidth}
 \scalebox{0.5}{\input{pd_bajd3m}}
 \caption{\label{fig:dens3m}Transition density $\Delta = 3/12$}
\end{subfigure}
\begin{subfigure}[c]{0.46\textwidth}
 \scalebox{0.5}{\input{pd_bajd2m}}
 \caption{\label{fig:dens2m}Transition density $\Delta = 2/12$}
\end{subfigure} 
\caption{\label{fig:dens}{\bf Comparison of density approximations. } This figure shows transition density approximations of the Basic Affine Jump-Diffusion (BAJD) solving the stochastic differential equation $dY_t=(\kappa \theta - \kappa  Y_t)\,dt+\sigma \sqrt{Y_t}dW_t +dL_t$. The parameters used are $\kappa \theta = 0.05, \kappa = 1, \sigma = 0.2, \lambda=1, \nu=0.05, y_0=0.05$. Panels \subref{fig:dens3m} and \subref{fig:dens2m} show the approximation for a time span of $\Delta = 3/12$, and $\Delta=2/12$, respectively. For both pictures, the exact density is obtained from Fourier inversion, while the density approximation FMS facilitates the approach from \citet{filipovicmayerhoferschneider13}.}
\end{figure}
Figure \ref{fig:dens} shows that the FMS density becomes negative close to zero, the more negative, the smaller $\Delta$, as Panels \ref{fig:dens3m} and \ref{fig:dens2m} indicate. Thus, in an application demanding positivity, such as derivatives pricing (\citet{aitsahalia99}), likelihood ratio tests, or MCMC sampling, it is imperative to use the PLR proposed in this paper, rather than the FMS projections without modifications.

\subsection{Distributions from options}\label{sec:optionpricing}
Deviations from positivity of probability distributions may have sizable economic ramifications. 
In this section, we explore PLR in the context of option pricing and the concept of no-arbitrage. The absence of arbitrage, or free lunch, implies that any payoff $X _{t+\Delta}$ has forward price $ F_t(X_{t+\Delta})$, where $ F _{t}$  is a linear operator that has a representation in terms of an expectation\footnote{We conduct our study in the forward market. Alternatively, in the spot market, we would assume zero interest rates for simplicity.}   
\begin{equation}\label{eq:price}
  F_t(X_{t+\Delta})=\mathbb E ^\QQ[X_{t+\Delta}|\mathcal F _t],
\end{equation} 
where $\QQ$ is a so-called forward measure. As before, we assume in this section that $\QQ$ is absolutely continuous with respect to the Lebesgue measure with density $q$ (and likewise $\PP$ with density $p$). Note that here $\PP$ represents an auxiliary possibly misspecified parametric probability measure approximating $\QQ$. 

A call option written on the payoff $X_{t+\Delta}$ with strike price $K$  pays $(X_{t+\Delta}-K)^+$, while for  puts the payoff is $(K-X_{t+\Delta})^+$. Given a set $\mathcal K:=\{K_1, \ldots, K_\Omega \}$ of strike prices along with bid-ask spreads, or for simplicity, mid quotes of option forward prices $P_{t,\Delta}(K_1),\ldots,P_{t,\Delta}(K_\omega),C_{t,\Delta}(K_{\omega+1}),\ldots, C_{t,\Delta}(K_{\Omega})   $ of out-of-the-money options (meaning that $K_1<\cdots <K_\omega\leq F_t(X_{t+\Delta})\leq K_{\omega+1}<\cdots<K_\Omega$), our goal is to find a density $q$ that is compatible with these prices. If strike prices were continuously quoted, $q$ could be computed from second derivatives of option prices \citep{breedenlitzenberger78}, but the presence of only a discrete observation $\mathcal K$ requires additional considerations to obtain $q$.

Option cross sections are informative financial instruments, as from \citet{carrmadan01}, for any twice-differentiable almost-everywhere function $f$, its forward-neutral expectation is a particular option portfolio,
\begin{equation}\label{eq:carrmadan}
 \mathbb E^\QQ[f(X_{t+\Delta})|\mathcal F_t]=\int _{0}^{F_t(X_{t+\Delta})}f''(K)P_{t,\Delta}(K)dK+\int _{F_t(X_{t+\Delta})}^{\infty}f''(K)C_{t,\Delta}(K)dK,
\end{equation}
so that, for instance, forward-neutral moments can be approximated from observed option cross sections. We will use formula  \eqref{eq:carrmadan} to compute \(\mu_i^\QQ:=\mathbb E^\QQ[ (\log X_{t+\Delta})^i|\mathcal F_t]\) and $\sigma ^\QQ:=\sqrt{\mu_2^\QQ-\mu_1^\QQ}$ to define
\begin{equation}
 Y:= \frac{\log \frac{X_{t+\Delta}}{F_t(X_{t+\Delta})}-\mu_{1}^\QQ}{\sigma ^\QQ}, \text{ and similarly } w:= \frac{\log \frac{K}{F_t(X_{t+\Delta})}-\mu_{1}^\QQ}{\sigma ^\QQ}.
\end{equation} 
From this change of variables, forward option prices can be written as
\begin{equation}
\begin{split}
 P_{t,\Delta}(K)&=\mathbb E ^{\QQ}[(K-X_{t+\Delta})^+]=e^ {\mu_1^\QQ} F_t(X_{t+\Delta})\int _{-\infty}^{\log w}(e^{\sigma ^\QQ w}-e ^{\sigma ^\QQ y})\tilde q(y)dy, \text{ and} \\
  C_{t,\Delta}(K)&=\mathbb E^{\QQ}[(X_{t+\Delta}-K)^+]=e ^{\mu_1^\QQ} F_t(X_{t+\Delta})\int _{\log w}^{\infty}(e ^{\sigma^\QQ y}-e^{\sigma^\QQ w})\tilde q(y) dy,
 \end{split}
\end{equation} 
where the density $\tilde q$ pertains to $Y$. From $\tilde q$, we can recover the original density $q$, the object of interest,  from  
\begin{equation*}
 q(x)=\frac{1}{x \, \sigma ^Q}\tilde q\brac{\frac{  \log \frac{x}{F_t(X_{t+\Delta})}-\mu_1^Q}{\sigma ^Q}}.
\end{equation*}
The above centering and scaling facilitate the object of interest
\begin{equation*}
 \xi _n=\frac{\tilde q}{p}\in L^{2}_\PP,
\end{equation*}
such that the options are priced, and  martingale and normalization restrictions, are satisfied. In terms of the auxiliary density $p$, the martingale and density normalization restrictions that should be satisfied under the forward measure read 
\begin{align*}
 1=e^ {\mu_1^Q} \int _{-\infty}^{\infty}e^{\sigma ^Q y} \xi_n(y)p(y)dy,\text{ and } 1=\int _{-\infty}^{\infty}\xi_n(y)p(y)dy, \text{ respectively}.
\end{align*}

In the context of fat-tailed forward-neutral $\tilde q$ densities, the choice of the auxiliary density $p$ becomes particularly important. \citet{lee04} investigates the tail of the distribution of forward-neutral log returns showing that it is exponential, and therefore ticker than the Gaussian tail. Thus, in the context of option pricing, the absolute continuity of $p$ with respect to $\tilde q$ puts strong restrictions on the choice of $p$, excluding in particular the normal distribution, as its tails are too thin. 

From these considerations, we solve the problem using the generalized hyperbolic distribution \citep{barndorff-nielsenhalgreen77} as auxiliary density $p$. The generalized hyperbolic distribution has four parameters. We choose the parameters to fit the first four moments of $Y$, such that $\mu^{\tilde  \QQ}_1=0,  \mu ^{\tilde  \QQ}_2=1, \mu^{\tilde  \QQ}_3, \mu ^{\tilde  \QQ}_4$. These moments can be computed either from a model under consideration, or using formula \eqref{eq:carrmadan} with observed option prices. Importantly, the generalized hyperbolic distribution features exponential tails, such that the approximation in $L^{2}_\PP$ is meaningful. In particular, for the parameterizations used below,  the tails of the generalized hyperbolic distribution satisfy the $L^{2}_\PP$ integrability condition.

Specializing program \eqref{eq:theprogram} for the task of extracting a density from option prices, we    
\begin{equation}\label{eq:theoptionprogram}
\begin{split}
 &\underset{  \xi_n \in M_n(D)}{\text{minimize }}   \norm{\xi_n}^2, \\
 &\text{ subject to } (\xi _n,1)=1, \,  (\xi _n,e^ {\mu_1^Q} e^{\sigma ^Q \cdot Y})=1, \\
  &1-\varepsilon\leq  \left (\xi_n , \frac{e^ {\mu_1^Q} F_t(X_{t+\Delta})(e^{\sigma ^Q \cdot w_j}-e ^{\sigma ^Q \cdot Y})^+}{P_\Delta(K_j)}\right )\leq 1+\varepsilon, \, j=1,\ldots ,  \omega \\
  &1-\varepsilon \leq \left (\xi_n , \frac{e^ {\mu_1^Q} F_t(X_{t+\Delta})(e^{\sigma ^Q \cdot Y}-e ^{\sigma ^Q \cdot w_{j}})^+}{C_\Delta(K_j)} \right )\leq 1+ \varepsilon, \,  j=\omega+1,\ldots , \Omega 
  \end{split}
\end{equation}
Tolerance $\varepsilon >0$ bounds the maximal relative option pricing error from above. \citet{luqu21} perform a similar task with polynomials by minimizing the squared pricing errors. Differently from their approach, we obtain a \emph{positive} density $q$, that furthermore satisfies the martingale restriction $(\xi _n,e^ {\mu_1^Q} e^{\sigma ^Q \cdot Y})=1$. Additionally, the generalized hyperbolic distribution features exponential tails, while \citet{luqu21} employ the normal distribution, for which the approximating expansion  diverges asymptotically.

We illustrate our approach using the \citet{duffiepansingleton00} double-jump model for   $\log X_t$, and its stochastic variance process $V_t$ with joint  dynamics
\begin{equation}\label{eq:dpsprocess}
 d\begin{pmatrix}
   \log X_t \\
   V_t
  \end{pmatrix}
  =
  \begin{pmatrix}
   r-\bar \lambda \bar \mu -\frac{1}{2}V_t \\
   \kappa _v(\bar v-V_t)
  \end{pmatrix}dt
+
\sqrt{V_t}
\begin{pmatrix}
 \rho & \sqrt{1-\rho^2} \\
 \sigma _v & 0
\end{pmatrix}
d\begin{pmatrix}
 W^{q}_{1,t} \\
 W^{q}_{2,t}
\end{pmatrix}
+dJ ^Q_t,
\end{equation}
where $W^{\QQ}_{1,t},W^{\QQ}_{2,t}$ are $\QQ$ Brownian motions, and $J ^\QQ$ is a pure jump process in $\mathbb R^2$ with constant mean arrival rate $\bar \lambda$, and $\bar \mu$ is the mean jump size in the $\log X_t$ direction under $\QQ$. 

For fixed $\Delta$, the characteristic function of $\log X_{t+\Delta}, V_{t+\Delta}|\log X_{t}, V_{t}$ is readily available in closed form, such that option prices can be obtained via the dampened transform as proposed by  \citet{lee04}. Using the estimated parameters from \citet[][TABLE I]{duffiepansingleton00}, we consider three scenarios for the value of $V_t$: low ($\sqrt{V_t}=0.02$), medium  ($\sqrt{V_t}=0.087$, the value quoted in the original paper), and high  ($\sqrt{V_t}=0.17$), and compute the resulting densities with positivity restriction, and without. For both the positive and the unrestricted PLR we find the smallest possible tolerance $\varepsilon$ that makes program \eqref{eq:theoptionprogram} feasible through trial and error. This process takes a fraction of a second.  

\begin{figure}
\centering
\begin{subfigure}[c]{0.46\textwidth}
 \scalebox{0.5}{\input{pd_dps_dens_low}}
 \caption{\label{fig:denslowvol}Transition density, low volatility}
\end{subfigure}
\begin{subfigure}[c]{0.46\textwidth}
 \scalebox{0.5}{\input{pd_dps_dens_med}}
 \caption{\label{fig:densmedvol}Transition density, medium volatility}
\end{subfigure} \\
\begin{subfigure}[c]{0.46\textwidth}
 \scalebox{0.5}{\input{pd_dps_dens_high}}
 \caption{\label{fig:denshighvol}Transition density, high volatility}
\end{subfigure}
\caption{\label{fig:dpsdens}{\bf Densities implied from options.} The figure shows the true density, the density solving program \eqref{eq:theoptionprogram} without (no restr.), and with (pos.) positivity restriction for the \citet{duffiepansingleton00} double-jump model estimated from options. Option data is generated for low ($\sqrt{V_t}=0.02$), medium  ($\sqrt{V_t}=0.087$), and high  ($\sqrt{V_t}=0.17$) volatility at valuation time. Both types of PLR are computed with polynomial degree $n=8$. In the horizontal axis we observe log returns $log \frac{X_{t+\Delta}}{X_t}$ of the underlying asset.}
\end{figure}
Figure \ref{fig:dpsdens} shows transition densities computed from Fourier inversion of model-implied option prices (true), without (no restr.), and with (pos.) positivity restriction imposed when solving program \eqref{eq:theoptionprogram}. Panels \ref{fig:denslowvol} and \ref{fig:denshighvol} show reasonably accurate $q$ densities for the low and high volatility scenarios, while the medium volatility in Panel \ref{fig:densmedvol} scenario is estimated too peaked. The reason for this is the auxiliary generalized hyperbolic distribution $p$, which despite matching the first four moments of the original distribution $\tilde q$, features a  peak that an order eight polynomial can not tilt sufficiently. Note that the unrestricted transition density achieves large negative values for negative returns close to zero.  

\begin{figure}
\centering
\begin{subfigure}[c]{0.46\textwidth}
 \scalebox{0.5}{\input{pd_dps_lr_low}}
 \caption{\label{fig:lrlowvol}$\xi _n$, low volatility}
\end{subfigure}
\begin{subfigure}[c]{0.46\textwidth}
 \scalebox{0.5}{\input{pd_dps_lr_med}}
 \caption{\label{fig:lrmedvol}$\xi _n$, medium volatility}
\end{subfigure} \\
\begin{subfigure}[c]{0.46\textwidth}
 \scalebox{0.5}{\input{pd_dps_lr_high}}
 \caption{\label{fig:lrhighvol}$\xi _n$, high volatility}
\end{subfigure}
\caption{\label{fig:dpslr}{\bf Polynomial tilts. } The figure shows the resulting polynomial $\xi _n$ from the solution of the program  \eqref{eq:theoptionprogram} without (no restr.), and with (pos.) positivity restriction for the \citet{duffiepansingleton00} double-jump model estimated from options. Option data is generated for low ($\sqrt{V_t}=0.02$), medium  ($\sqrt{V_t}=0.087$), and high  ($\sqrt{V_t}=0.17$) volatility at valuation time. Both types of PLR are computed with polynomial degree $n=8$. In the horizontal axis we observe log returns $log \frac{X_{t+\Delta}}{X_t}$ of the underlying asset.}
\end{figure}

Figure \ref{fig:dpslr} illustrates the stark differences between the unrestricted PLR and the positive one induced by the positivity condition. It also shows that the tilting is moderate in magnitude, thanks to the minimum-norm objective function. The unrestricted $\xi _n$, i.e., the one with no positivity restriction, takes negative values for all scenarios but the high-volatility one and presents an excessive oscillatory behavior when volatility is low or medium. 

\begin{figure}
\centering
\begin{subfigure}[c]{0.46\textwidth}
 \scalebox{0.5}{\input{pd_dps_iv_low}}
 \caption{\label{fig:ivlowvol}Implied volatility, low volatility}
\end{subfigure}
\begin{subfigure}[c]{0.46\textwidth}
 \scalebox{0.5}{\input{pd_dps_iv_med}}
 \caption{\label{fig:ivmedvol}Implied volatility, medium volatility}
\end{subfigure} \\
\begin{subfigure}[c]{0.46\textwidth}
 \scalebox{0.5}{\input{pd_dps_iv_high}}
 \caption{\label{fig:ivhighvol}Implied volatility, high volatility} 
\end{subfigure}
\caption{\label{fig:dpsiv}{\bf Comparison of implied volatilities. }The figure shows Black-Scholes implied volatilities of the option prices generated by the true distribution (true), the density solving program \eqref{eq:theoptionprogram} without (no restr.), and with (pos.) positivity restriction for the \citet{duffiepansingleton00} double-jump model estimated from options as a function of log moneyness. Option data is generated for low ($\sqrt{V_t}=0.02$), medium  ($\sqrt{V_t}=0.087$), and high  ($\sqrt{V_t}=0.17$) volatility at valuation time. Both types of PLR are computed with polynomial degree $n=8$. In the horizontal axis we observe log moneyness $log \frac{K}{X_t}$.}
\end{figure}
Figure \ref{fig:dpsiv} shows Black-Scholes implied volatilities from the true model, as well as the polynomial approximations. From Panel \ref{fig:ivlowvol} it becomes apparent that the low volatility scenario is the most challenging. In this scenario, OTM call options have low prices dominating the relative pricing error on both polynomial likelihood approximations what forces them to sacrifice part of their ability to price OTM puts. The consequences for the implied volatility curve is that both models underprice puts (negative log-moneyness in the picture) with the unconstrained model having smaller bias and an intermediary region of log-moneyness where it overprices puts. The positive model trades precision to avoid the excessive oscillatory behavior produced by the unconstrained PLR on the implied volatility curve, which comes from the oscillatory behavior of its estimated transition density noted in Figure \ref{fig:dpsdens}. This excessive oscillatory behavior of implied volatilities also appears in the medium volatility scenario for the unconstrained model, contrasting with well-behaved implied volatilities under positivity constrain. Finally, note that both approximations work quite well in the high volatility scenario.

\section{Conclusion}\label{sec:conclusion}
We develop projections of likelihood ratios onto polynomials that preserve positivity. We term them positive polynomial likelihood ratio (PLR). PLR can accommodate shape restrictions, are fast and robust to compute as solutions to conic programs, and come with asymptotic theory for their use with sample moments.

We illustrate  PLR with two applications. The first is an approximation of the unknown transition density of a jump-diffusion process. The second construct a density from the observation of option prices only. In both cases, the virtue of positivity of the PLR becomes evident through the positivity of the resulting distributions, that without the restriction, becomes negative and with excessive oscillatory behavior. 

PLR lend themselves to nonparametric structurally constrained settings, such as for instance  \citet{christensen17} discount factor decomposition. We leave this, and related applications to future research.

\bibliographystyle{plainnat}  
\bibliography{master} 

\begin{appendix}
\appendixpage

\section{Positive polynomials}\label{sec:pospol}
\counterwithin{theorem}{section}
\renewcommand{\thetheorem}{A\arabic{theorem}}
In this section, we review results in the literature about positive polynomials.
In the univariate case, we have 
\begin{proposition}[\citet{schmuedgen17}]\label{prop:pospol}
 For  any positive integer $n$,
\begin{enumerate}
 \item $D=\reals$: $\Omega_{n}:=M_{2n}=\brat{f^2(t)+g^2(t):f,g\in \reals [t]_n}$,
 \item $D=\reals _+$: $M_{2n}=\brat{f(t)+t g(t) :f\in \Omega _n, g \in \Omega _{n-1}}$,
 \item $D=\reals _+$: $M_{2n+1}=\brat{f(t)+t g(t) :f,g\in \Omega _n}$,
 \item $D=[a,b]$: $M_{2n}=\brat{f(t)+(b-t)(t-a)g(t):f\in \Omega _n, g \in \Omega _{n-1}}$,
 \item $D=[a,b]$: $M_{2n+1}=\brat{(b-t)f(t)+(t-a)g(t):f,g\in \Omega _n}$.
\end{enumerate}
\end{proposition}
The set of positive polynomials on any other (continuous) state space can be extracted from Proposition \ref{prop:pospol} from a change of variables. For instance, $D=[a,\infty)$ can be obtained from parameterization (2) above through the change of variables $p(t-a)$ for $p\in M_{2n}$ on $D=\reals _+$. 

In the multivariate case, nonnegative polynomials exist that are not s.o.s.. Since we merely want to assure non-negativity,  and a s.o.s. polynomial is certainly non-negative, it is sufficient for our purpose to work with s.o.s. polynomials. Any such polynomial has a representation as a quadratic form (the proof is in \citet{schmuedgen17} for Proposition 13.2).
\begin{proposition}\label{prop:posmultivar}
 A polynomial $\xi_{2n}\in \reals[t]_{2n}$ is s.o.s. if and only if
  $\xi_{2n} =\tau_n ^{\T}\, \bm V \, \tau_n, \text{ with } \bm V\succeq 0. $
\end{proposition}
Note that we do not make a distinction between different supports in the multivariate case. 
Statements   \ref{prop:pospol} and \ref{prop:posmultivar}   lead to our parameterization \ref{def:MD} of the cone of positive polynomials.

\section{Proofs}\label{sec:proofs}

\noindent \textit{Proposition \ref{prop:unique}}
\begin{proof}
 For the minimization, it is convenient to use one-half the squared norm as an objective function rather than the norm itself. This does not change the result, since the norm is non-negative, and we can write 
\begin{equation*}
 \norm{\xi_n}^{2}/2=\frac{1}{2}(\xi_n, \xi_n)=\frac{1}{2}\int _{D}\bm x ^{\T}\bm \tau_n(t) \bm \tau_n^{\T}(t)\bm x \, d\PP(t)=\frac{1}{2}\bm x ^{\T}\bm H_n\bm x .
\end{equation*}
Together with the constraints in coordinate form and Proposition \ref{prop:posmultivar}, this yields primal \eqref{eq:mixedconicprimal} that can be solved as a mixed conic semidefinite program.\footnote{We use the \href{https://www.mosek.com}{Mosek} optimizer to solve the program in practice.} 
 Since with Assumption \ref{ass:card} the objective function is strictly convex, and the constraint set is an intersection of closed convex sets,  the solution 
 with Assumption \ref{ass:feasible} is unique in $\bm x$, and strong duality obtains. With strong duality at hand, we next consider the dual form to \eqref{eq:mixedconicprimal}. 
Furthermore, the cone of symmetric positive semidefinite matrix is self-dual. From these observations,   we can write the Lagrangian of system \eqref{eq:mixedconicprimal} as 
\begin{equation}\label{eq:mylagrangian}
\begin{aligned}
 \mathcal L(\bm x,\bm \eta, \bm \nu, \bm \varepsilon ,\bm \Lambda _V, \bm \Lambda _W )&:=\frac{1}{2}\bm x ^{\T}\bm H_n\bm x -\bm \eta ^{\T}(\bm F_n \bm x-\bm c)-\bm \nu ^{\T}(\bm x-T(\bm V, \bm W))\\
 &+\bm \varepsilon ^{\T}(\bm G_n \bm x-\bm d)-\tr (\bm \Lambda_V \bm V)-\tr (\bm \Lambda _W \bm W)
 \end{aligned}
\end{equation} 
for $\bm \eta\in \reals ^{m}, \bm \varepsilon \in \reals _+^l, \bm \nu \in \reals ^{N+1}$, as well as p.s.d. matrices $\bm \lambda _V, \bm \lambda _W$ of the same dimension as $\bm V$ and $\bm W$, respectively. 
From the first-order condition (on $\bm x$) we can then deduce
 $\bm x_0 ^{\T}\bm H_n - \bm \eta ^{\T}_0\bm F_n -\bm \nu ^{\T}_0+ \bm \varepsilon  ^{\T}_0\bm G_n=0$,
and with $\bm H_n$ invertible from Assumption \ref{ass:card}, any optimal solution must satisfy \eqref{eq:myfoc}. The Lagrangian is linear-quadratic  in $\bm V$ and $\bm W$, and taking matrix derivatives, we therefore  have the conditions 
 $\nu _0 \bm V_0+\cdots + \nu _N \bm V_N=\bm \Lambda _V\succeq 0$, and analogously  $\nu _0 \bm W_0+\cdots +  \nu _N \bm W_N=\bm \Lambda _W\succeq 0$.  Plugging these relations into the Lagrangian yields
 \begin{align*}
  g(\bm \eta, \bm \nu, \bm \varepsilon)&=\frac{1}{2}(\bm \eta ^{\T} \bm F_n  + \bm \nu^{\T} -\bm \varepsilon ^{\T}\bm G _n)\bm H_n ^{-1}(\bm F_n ^{\T}\bm \eta + \bm \nu -\bm G _n^{\T}\bm \varepsilon) \\
  &-\bm \eta ^{\T}(\bm F_n \bm H_n ^{-1}(\bm F_n ^{\T}\bm \eta + \bm \nu -\bm G _n^{\T}\bm \varepsilon)-\bm c) \\
  &-\bm \nu ^{\T}\bm H_n ^{-1}(\bm F_n ^{\T}\bm \eta + \bm \nu -\bm G _n^{\T}\bm \varepsilon) \\
  &+\bm \varepsilon ^{\T}(\bm G_n \bm H_n ^{-1}(\bm F_n ^{\T}\bm \eta + \bm \nu -\bm G _n^{\T}\bm \varepsilon)-\bm d).
 \end{align*}
The objective function in \eqref{eq:mixedconicprimal}  can be obtained by completing the square.
\end{proof}

\noindent \textit{Proposition \ref{pro:consistency}}  

\begin{proof} 
In this section, we denote by $\|\cdot \|$ the Euclidean norm for lighter notation (rather than the $L^{2}_\PP$ norm). 
 For convenience, define $\RR_k(\bm x):=\bm x^{\T}\hat{\bm H}_{n,k}\bm x$, $\RR_0(\bm x):=\bm x^{\T}\bm H_n \bm x$ and
denote by $\bm x _0$ the solution of \eqref{eq:mixedconicprimal}. For a set $A\subseteq \mathbb{R}^{N+1}$ and $\bm x\in\mathbb{R}^{N+1}$, define the distance $d(\bm x,A):=\inf_{\bm y\in A} \|\bm x-\bm y\|$ and the $\epsilon$-enlargement of $A$ as $A^{\epsilon}:=\{x\in\mathbb{R}^{N+1}:d(\bm x,A)\leq \epsilon\}$. Let $\bm v_0$ denote a  vector collecting all the vectorized population quantities ($\vecme \bm H_{n},\vecme \bm F_{n},\ldots$ ) to be estimated from the data and  $\widehat{\bm v}_k$ its respective estimator. Finally For $\epsilon, \rho>0$, define the event $\mathcal{A}_k:=\{\|\widehat{\bm v}_k-\bm v_0\|\leq \epsilon, \max(\|\bm x_0\|,\|\widehat{\bm x}_k\|) \leq\rho, \sup_{\bm \|x\|\leq \rho}|\RR_k(\bm x) - \RR_0(\bm x)|\leq \epsilon\}$. 

We claim that the event $\mathcal{A}_k$ implies
\begin{enumerate}
\item[(i)] $\tilde K \subseteq\tilde K _k^{\rho\epsilon}$; 
\item[(ii)] $\Delta_\epsilon:=\sup_{\bm x\in \widetilde{K}_k^{\rho\epsilon}, \bm y\in \widetilde{K}_k}|\RR_k(\bm x) - \RR_k(\bm y)|\leq 2\rho^2(\|\bm H_n\| +\epsilon)\epsilon$ for every $\delta>0$;
\item[(iii)] $\RR_0(\widehat{\bm x}_k) - \RR_0(\bm x_0)\leq 2\left[1 +\rho^2(\|\bm H_n\| +\epsilon)\right]\epsilon$.
\end{enumerate}
Indeed, for $\bm x\in\widetilde{K}$ we have $d(\bm x,\widetilde{K}_k)\leq \|\widehat{\bm v}_k-\bm v_0\|\|\bm x\|\leq \rho\epsilon $  which shows (i). For (ii), we have for $\bm x\in \widetilde{K}_k^{\rho\epsilon}$ and $\bm y\in \widetilde{K}_k$
\[
|\RR_k(\bm x) - \RR_k(\bm y)| = |(\bm x + \bm y)^{\T} \widehat{\bm H}_{n,k}(\bm x - \bm y)|\leq \|\bm x+\bm y\|_2\|\widehat{\bm H}_{n,k}\|_2\|\bm x-\bm y\|_2\leq 2\rho^2\|\widehat{\bm H}_{n,k}\|\epsilon,
\]
and $\|\widehat{\bm H}_{n,k}\|\leq \|\bm H_n\| + \|\widehat{\bm H}_{n,k} -\bm H_n\|\leq \|\bm H_n\| +\|\widehat{\bm v}_k-\bm v_0\|\leq \|\bm H_n\| +\epsilon$. Finally, for (iii)  we have
\begin{align*}
    \RR_0(\widehat{\bm x}_k) - \RR_0(\bm x_0) &= \RR_0(\widehat{\bm x}_k) - \RR_k(\bm x_0)  + \RR_k(\bm x_0) - \RR_0(\bm x_0)\\
    &\leq \RR_0(\widehat{\bm x}_k) - \RR_k(\widehat{\bm x}_k)  + \RR_k(\bm x_0) - \RR_0(\bm x_0) +\Delta_\epsilon\\
    &\leq 2\left[\sup_{\bm x}|\RR_k(\bm x) - \RR_0(\bm x)| +\rho^2(\|\bm H_n\| +\epsilon)\epsilon\right]\\
    &\leq 2\left[1 +\rho^2(\|\bm H_n\| +\epsilon)\right]\epsilon,
\end{align*}
where the first equality follows because $\RR_k(\bm x_0)\geq \RR_k(\bm x^*)$ for some $\bm x^*\in \widetilde{K}_k^{\rho\epsilon}\supseteq \widetilde{K}$ by (i); and $\RR_k(\bm x^*) \geq  \RR_k(\widetilde{\bm x}) - \Delta_\delta$ for some $\widetilde{\bm x}\in \widetilde{K}_k$ by (ii). Finally $\RR_k(\widetilde{\bm x})\geq \RR_k(\widehat{\bm x}_k)$ by the optimality of $\widehat{\bm x}_k$ for $\RR_k$ restricted to $\tilde K _k$.

Fix an arbitrary $\eta>0$. Due to the uniqueness of the solution (Proposition \ref{prop:unique}) and continuity of the objective function, the event $\|\widehat{\bm x}_k-\bm x_0\|> \eta$ implies that $\RR_0(\widehat{\bm x}_k) - \RR_0(\bm x_0)> \delta$ for some $\delta>0$. For any given $\rho$ we always take $\epsilon$ small enough in (iii) to have $\RR_0(\widehat{\bm x}_k) - \RR_0(\bm x_0)\leq \delta$ on $\mathcal{A}_k$. Therefore 
\begin{align*}  
\PP(\|\widehat{\bm x}_k-\bm x_0\|_2> \eta)&\leq  \PP(\RR_0(\widehat{\bm x}_k) - \RR_0(\bm x_0)>\delta)\\
&\leq \PP(\{\RR_0(\widehat{\bm x}_k) - \RR_0(\bm x_0)> \delta\}\cap \mathcal{A}_k) + \PP(\mathcal{A}_k^c)\\
&\leq \PP(\|\widehat{\bm v}_k-\bm v_0\|> \epsilon) + \PP(\max(\|\bm x_0\|,\|\widehat{\bm x}_k\|) >\rho ) \\
&\qquad +\PP(\sup_{\bm \|x\|\leq \rho}|\RR_k(\bm x) - \RR_0(\bm x)|> \epsilon).
\end{align*}

The first term in the last expression can be made arbitrarily small by taking $k$ large enough by Assumption \ref{ass:rate}. Also, since $\RR_0$ is strictly convex by Assumption \ref{ass:card} and $\widehat{\bm H}_{n,k}\to \bm H_n$ in probability as $k\to\infty$ for fixed $n$ by Assumption \ref{ass:rate}, have that $\{\RR_k\}_k$ is strongly convex with  high probability (for $k$ large enough). Therefore there exist a closed all $B\subseteq\mathbb{R}^{N+1}$ of radius $\rho>0$ and $N_0\in\mathbb{N}$  such that $\bm x_0,\widehat{\bm x}_k \in B$ for all $k\geq N_0$  with high probability. Thus the second term can be made arbitrarily small, by taking the radius $\rho$ large enough. Finally, $\{\RR_k\}$ is a sequence of strictly convex functions (with high probability) so it converges to $\RR_0$ in probability uniformly on a compact set under Assumption \ref{ass:rate}. 
That concludes the consistency proof.
\end{proof}  

\noindent \textit{Proposition \ref{pro:rate}}   

\begin{proof} In this section, we denote by $\|\cdot \|$ the Euclidean norm for lighter notation (rather than the $L^{2}_\PP$ norm). Fix $\epsilon>0$ and pick a compact $\varXi \subseteq \Theta_*\cup\widehat{\Theta}$ such that $\|\widehat{\bm s} - \bm s\|\leq C_\varXi  r_k$ with probability at least $1-\epsilon$. First we note that $\Psi$ uniformly convergent in probability over $\varXi \subseteq \Theta$, because
\begin{equation*}
 \sup_{\bm \theta\in \varXi }\|\widehat{\Psi}(\bm\theta) - {\Psi}(\bm\theta)\|\leq C_1\|\widehat{\bm s} - \bm s\|\sup_{\bm \theta\in K}\|\bm \theta\|\leq C_2\|\widehat{\bm s} - \bm s\|\leq C_\epsilon r_k.
\end{equation*}
where $C_1$ is a constant depending on the dimensions $N, \ell, m$; $C_2:=C_1\text{diam}(\varXi )$ and $C_\epsilon:= C_2C_\varXi $. From here on and below are conditioning on the event that $\{\|\widehat{\bm s} - \bm s\|\leq C_\varXi  r_k\}$. For $\bm\theta\in\Theta_*$ we have $\|\widehat{\Psi}(\bm\theta)\| =  \|\widehat{\Psi}(\bm\theta) - {\Psi}(\bm\theta)\|\leq Cr_k$ thus $\Theta_*\subseteq \widehat{\Theta}$ provided that $\delta\geq C_\epsilon r_k$. 

Similarly, for $\bm\theta \in\widehat{\Theta}$ we have $\|\Psi(\bm\theta)\| \leq   \|\widehat{\Psi}(\bm\theta) - {\Psi}(\bm\theta)\| + \|\Psi(\bm\theta)\| \leq C_\epsilon r_k + \delta$. Now, fix $\eta>0$ and let $\bm\theta \in\widehat{\Theta}\setminus\Theta_*^\eta$ and $\bm\theta_*$ be the (unique) point closest to $\bm\theta$ that belongs to $\Theta_*$. Finally,  let $\widetilde{\bm\theta}$ be the point in the middle of the line segment between $\bm\theta$ and  $\bm\theta_*$. Since $\widetilde{
\bm\theta}\notin \Theta_*$, we have that $\|\Psi(\widetilde{\bm \theta})\|>0$ and  $\nabla\|\Psi(\widetilde{\bm\theta}) \|^2\geq C_0$ for some constant $C_0>0$. The last condition holds because 
$\|\Psi(\bm\theta)\|^2$ is a differentiable convex non-negative function with minimum when $\Psi(\bm\theta)=0$ (equivalently when $\bm\theta \in\Theta_*$). Moreover $\Psi(\bm\theta)=0$ if and only if $\nabla \|\Psi(\bm\theta)\|^2=\bm 0$. Therefore, $\bm\theta \in\widehat{\Theta}\setminus\Theta_*$ implies both $\|\Psi(\bm \theta)\|>0$ and $\nabla\|\Psi(\bm\theta)\|\neq \bm 0$ as stated. Also, due to the convexity of $\|\Psi(\bm\theta)\|^2$ we have $\|\Psi(\bm \theta)\|^2 \geq \|\Psi(\widetilde{\bm \theta})\|^2 + \nabla \|\Psi(\widetilde{\bm\theta})\|^2(\widetilde{\bm \theta} - \bm\theta)> 0 + C_0\eta/2$. Take $\eta\to 0$ and $k\to\infty$ and eventually (for large k) $C_0\eta/2\leq 1$ and we have $\|\Psi(\bm \theta)\| >\sqrt{C_0\delta/2}\geq C_0\eta/2$. 

Set $\delta = C r_k$ for $C\geq C_\epsilon$ and $\eta = (C/C_0)r_k$ to conclude that both inequalities imply that if $\bm\theta \in\widehat{\Theta}$ then $\bm\theta\in \Theta_*^{Mr_k}$ for large $k$ and some positive constant $M:=M_\epsilon$ as $k\to \infty$. We can choose $M_\epsilon$ such that $\widehat{\Theta}\subseteq \Theta_*^{M_\epsilon r_k}$ for every $k\in \mathbb{N}$. Therefore, $d(\widehat{\Theta},\Theta_*)=\inf\{\omega\geq 0:\Theta_*\subseteq \widehat{\Theta}^\omega, \Theta_*\subseteq \widehat{\Theta}^\omega\}\leq M_\epsilon r_k $ with probability at least $1-\epsilon$.

A careful review of the proof of Proposition~\ref{pro:rate} reveals  Corollary \ref{cor:coverage}, which can be useful for inferring coverage of $\Theta_*$ by set estimator $\widehat{\Theta}$.

\end{proof}

\noindent \textit{Theorem \ref{th:L2}}

\begin{proof}
In this application, the subspace $K$ is itself generated by monomials. This allows a direct sum decomposition 
\begin{equation}
 \xi _n^{(m)}=\xi ^{\star(m)}_n\oplus \xi ^{\circ(m)}_n,
\end{equation}
where $\xi ^{\star (m)}_n$ is the minimum-norm polynomial projection from \eqref{eq:fms}, and $\xi ^{\circ (m)}_n$ is the minimum-norm polynomial that lifts $\xi_n^{(m)}$ into the cone of pointwise positive polynomials. 
From Assumption \ref{ass:feasible}, the solution is feasible. From Assumption \ref{ass:basis}, $P_{\PP,n}\oplus P_{\PP,n}^{\perp}=L^{2}_\PP$. We can therefore write
 $\frac{d\QQ}{d\PP}=\xi^{\star(m)}_n+ \xi ^{\circ(m)}_n+ \epsilon$
with $\xi^{\star(m)}_n\in K, \xi ^{\circ(m)}_n\in K^{\perp}, \xi_n^{\star(m)}+\xi ^{\circ(m)}_n\in M_n(D)\cap K \subseteq P_{\PP,n}$, and $\epsilon \in P_{\PP,n}^{\perp}$. For each $m$,   $\xi ^{\star(m)}_n$ solves the standard Hilbert minimum-norm problem (cf. Section \ref{sec:rkhspzn}), and from Assumption \ref{ass:basis} it converges in $L^{2}_\PP$ to $d\QQ/\PP$ \citep{filipovicmayerhoferschneider13}. From this, we can write
 \begin{equation*}
  0=\lim _{m\to \infty}\norm{\frac{d\QQ}{d\PP}-\xi ^{\star(m)}_n}^2=\lim _{m\to \infty}\norm{\xi ^{\circ(m)}_n+\epsilon}^2=\lim _{m\to \infty}\norm{\xi ^{\circ(m)}_n}^2+\lim _{m\to \infty}\norm{\epsilon}^2.
 \end{equation*} 
 From the non-negativity of the norms  $\lim _{m\to \infty}\norm{\xi_n ^{\circ(m)}}=0$ and $\lim _{m\to \infty}\norm{\epsilon }=0$. Therefore $
  \lim _{m\to \infty}\norm{\frac{d\QQ}{d\PP}-\xi ^{(m)}_n}=\lim _{m\to \infty}\norm{\epsilon}=0$.
\end{proof}

 \end{appendix} 
\end{document}